\documentclass[a4paper,twoside,10pt]{article}

\usepackage{amsmath,amsthm,amssymb,amscd,xspace}
\DeclareSymbolFont{rsfs}{OMS}{rsfs}{m}{n}
\DeclareSymbolFontAlphabet{\rsfs}{rsfs}
\let\mathcal=\rsfs

\usepackage[all,ps]{xy} 

\newcommand\leftidx[3]{%
  {\vphantom{#2}}#1#2#3%
}

\usepackage[frenchb,english]{babel}
 
\newtheorem{thm}{Theorem}[section]
 
 \newtheorem {lemma}            [thm]{Lemma}
 \newtheorem {cor}              [thm]{Corollary}
 \newtheorem*{thm*}                  {Theorem}
 \newtheorem {prop*}                 {Proposition}
 \newtheorem*{lemma*}                {Lemma}
 \newtheorem*{cor*}                  {Corollary}
\theoremstyle{definition}

 \newtheorem*{def*}                  {Definition}
 \newtheorem*{ex*}                   {Example}
 \newtheorem*{xca*}                  {Exercise}
 \newtheorem*{xcas*}                 {Exercises}
 \newtheorem*{ax*}                   {Axiom}
\theoremstyle{remark}
 \newtheorem {rmk}              [thm]{Remark}

 \newtheorem*{rmk*}                  {Remark}
 \newtheorem*{note*}                 {Note}
 \newtheorem*{add*}                  {Addendum}

\newcommand{\Op}[1]{#1^{\scriptscriptstyle\mathrm{op}}}
\newcommand{\lidx}[2]{\leftidx{^{#1}}{#2}{\relax}}
\newcommand{\dummy}[1][]{{\_}_{\mathchoice%
 {}{}{\raisebox{-.28ex}{$\scriptscriptstyle #1$}}%
 {\raisebox{-.4ex}{$\scriptscriptstyle #1$}}}}
\newcommand{\tmap}[5]{#1\colon #2\to #3\colon #4\mapsto #5}
\newcommand{\dmap}[5]{\begin{array}[c]{@{}l@{\ }c@{\ }l@{\ }c@{\,}}%
                       #1\colon & #2 & \longrightarrow & #3 \\
                                & #4 & \longmapsto     & #5
                      \end{array}}
\newcommand{\map}[5]{%
 \if@display
  \dmap{#1}{#2}{#3}{#4}{#5}
 \else
  \tmap{#1}{#2}{#3}{#4}{#5}
 \fi}
\newcommand{\C}            {\mathbb{C}}
\newcommand{\tunit}        {\mathbb{I}}
\newcommand{\Id}           {\mathrm{Id}}
\newcommand{\aut}          {\mathrm{Aut}}
\newcommand{\End}          {\mathrm{End}}

\newcommand{\Hom}          {\mathrm{Hom}}
\newcommand{\HomCheck}     {\Hom\check{\,}}
\newcommand{\EndCheck}     {\End\check{\,}}
\newcommand{\Vect}         {\mathcal{V}\!\mbox{\itshape ect}}

\newcommand{\Com}          {\mathcal{C}\!\mbox{\itshape om}}
\newcommand{\Mod}          {\mathcal{M}\!\mbox{\itshape od}}
\newcommand{\Tr}           {\mathrm{tr}}
\newcommand{\qtr}          {\mathrm{qtr}}
\newcommand{\krn}          {\mathrm{k}}
\newcommand{\cok}          {\mathrm{c}}
\newcommand{\Krn}          {\mathrm{K}}
\newcommand{\Cok}          {\mathrm{C}}
\newcommand{\pk}           {{\mathrm{pk}}}
\newcommand{\grande}       {{\mathrm{tot}}}
\newcommand{\alphath}      {{$\alpha$\nobreakdash-\hspace{0pt}th}\xspace}
\newcommand{\adjunction}[2]{\!\!\xymatrix@C=1pc{#1\ar@{-|}[r]& #2}\!\!}

\newenvironment{sketch}
{\begin{proof}[Proof \textup{(\/}sketch\textup{)}]}
{\end{proof}}

\title{Tannaka Reconstruction for\\
       Crossed Hopf Group Algebras}
\author{Marco Zunino}
\date{\today}
\begin{document}

\maketitle

\begin{abstract}
 \noindent We provide an analog of Tannaka Theory for Hopf algebras in the
 context of crossed Hopf group coalgebras introduced by Turaev.
 Following Street and our previous work on the quantum double of
 crossed structures, we give a construction, via Tannaka Theory, 
 of the quantum double of crossed Hopf group algebras (not 
 necessarily of finite type).
\end{abstract}

\setcounter{tocdepth}{1}
\tableofcontents

\section*{Introduction}
Turaev~\cite{Tur-pi,Tur-CPC} generalized Reshetikhin-Turaev~\cite{RT1} 
invariants and the notion of \textsc{tqft} to the case of 
$3$\nobreakdash-\hspace{0pt}manifolds endowed with a homotopy 
class of maps to $K(G,1)$, where $G$ is a discrete group
(see also L{\^e} and Turaev~\cite{LeTur} and Virelizier~\cite{Virelizi2}).
 One of the key points in~\cite{Tur-CPC} is the notion of 
\textit{crossed Hopf $G$\nobreakdash-\hspace{0pt}coalgebra}. 
In the same way as one can use categories of representations 
of modular Hopf algebras to construct Reshetikhin-Turaev 
invariants of $3$\nobreakdash-\hspace{0pt}manifolds, 
one can use categories of representations of modular 
crossed Hopf $G$\nobreakdash-\hspace{0pt}coalgebras  
to construct homotopy invariants of maps from 
$3$\nobreakdash-\hspace{0pt}manifolds to 
the Eilenberg-Mac~Lane space $K(G,1)$. 
Similarly, $G$\nobreakdash-\hspace{0pt}coalgebras
are used by Virelizier~\cite{Virelizi2} to construct 
homotopy invariants which generalize Hennings invariant of
$3$\nobreakdash-\hspace{0pt}manifolds.

Roughly speaking, a crossed Hopf $G$\nobreakdash-\hspace{0pt}coalgebra $H$ 
is a family $\{H_{\alpha}\}_{\alpha\in G}$ of algebras endowed with 
a \textit{comultiplication}
$\Delta_{\alpha,\beta}\colon H_{\alpha\beta}\to H_{\alpha}\otimes H_{\beta}$, 
a \textit{counit} $\varepsilon\colon\Bbbk\to H_{1}$ 
(where $1$ is the unit of $G\/$), 
and an \textit{antipode}
$s_{\alpha}\colon H_{\alpha}\to H_{\alpha^{-1}}$. 
It is also required that $H$ is endowed with a family 
of algebra isomorphisms 
$\varphi_{\beta}\colon H_{\alpha}\to H_{\beta\alpha\beta^{-1}}$, 
the \textit{conjugation,} compatible with the above structures 
and such that $\varphi_{\beta\gamma}=\varphi_{\beta}\circ\varphi_{\gamma}$.
If $G=1$, then we recover the usual definition of a Hopf algebra. 
A \textit{universal $R$\nobreakdash-\hspace{0pt}matrix}
and a \textit{twist} for a $H$ are, respectively, families
$R=\bigl\{\xi_{(\alpha)}\otimes\zeta_{(\beta)}=R_{\alpha,\beta}\in 
H_{\alpha}\otimes H_{\beta}\bigr\}_{\alpha,\beta\in G}$ and
$\theta=\{\theta_{\alpha}\in H_{\alpha}\}_{\alpha\in G}$ 
satisfying axioms that explicitly involve the conjugation.
Properties of $G$\nobreakdash-\hspace{0pt}coalgebra are studied
in~\cite{Virelizi}.

In this article, following the constructions given, 
in the case of a Hopf algebra, by Joyal and Street~\cite{JS-Tannaka},
we consider which conditions a functor $F$ from a tensor
category $\mathcal{C}$ to the category of complex vector spaces
should satisfy to obtain a crossed group Hopf algebra $H$ such that
$\mathcal{C}$ will be isomorphic to the category of $H$-comodules.
In doing that, we extend Tannaka Theory to the crossed case.
By using our generalization~\cite{Zunino-1}
of Joyal and Street's center construction~\cite{JS}, we provide a
construction of the quantum co-double of a crossed Hopf group coalgebra 
similar to the construction given by Street in the case of a Hopf 
algebra~\cite{Street-double}. More in detail, starting from a crossed
Hopf group coalgebra $H$, we construct a coquasitriangular crossed
Hopf group coalgebra $D'(H)$ such that the category of comodules of
$D'(H)$ is equivalent to the center of the category of comodules of
$D(H)$. In particular, if $D'(H)$ is of finite type (i.e.\@ if every
$D'(H)_{\alpha}$ is finite-dimensional), then it is isomorphic to the
dual of $D(H^{\ast})$, where $H^{\ast}$ is the dual of $H$ and
$D(H^{\ast})$ is the generalization of Drinfeld quantum
double~\cite{Drn} introduced in~\cite{Zunino-2}.

We should obseve that, once the correct axioms are fixed, then the
proofs follow in a relatively simple way by generalizing the standard theory.
However, to determinate the precise conditions for a crossed category to be
Tannakian is a necessary step preliminary to other projects like, for instance, 
the quantization of the crossed structures (which is discussed in an article 
in preparation).

\section{Crossed Group Categories}
We recall some basic definitions about crossed structures,
in particular the definition of a crossed group category.\bigskip

\subsection{Basic Definitions}
A \textit{Crossed $G$\nobreakdash-\hspace{0pt}category}~\cite{Tur-CPC}
$\mathcal{C}$ is given by the following data:
\begin{itemize}
 \item a tensor category $\mathcal{C}\/$;
 \item a family of subcategories $\{\mathcal{C}_{\alpha}\}_{\alpha\in G}$
       such that $\mathcal{C}$ is disjoint union of this 
       family and that
       $U\otimes V\in\mathcal{C}_{\alpha\beta}$, for any
       $\alpha,\beta\in G$, 
       $U\in\mathcal{C}_{\alpha}$, and 
         $V\in\mathcal{C}_{\beta}\/$;
 \item a group homomorphism
       $\map{\Phi}{G}{\aut(\mathcal{C})}{\beta}{\Phi_{\beta}}$,
       the \textit{conjugation,} (where $\aut(\mathcal{C})$ is the group 
       of invertible strict 
       tensor functors from $\mathcal{C}$ to itself) such that 
       $\Phi_{\beta}(\mathcal{C}_{\alpha})=\mathcal{C}_{\beta\alpha\beta^{-1}}$
       for any $\alpha,\beta\in G$.
\end{itemize}
Given $\alpha\in G$, the subcategory $\mathcal{C}_{\alpha}$ is said the 
\textit{\alphath component} of $\mathcal{C}$ while the functors $\Phi_{\beta}$ 
are said \textit{conjugation isomorphisms.} Given $\beta\in G$ and an object 
$V\in\mathcal{C}_{\beta}$, the functor $\Phi_{\beta}$ is denoted by
$\lidx{V}{(\cdot)}$, 
as in~\cite{Tur-CPC}, or even by $\lidx{\beta}{\,(\cdot)}$.
The crossed $G$\nobreakdash-\hspace{0pt}category $\mathcal{C}$ is 
\textit{strict} when it is strict 
as a tensor category. When $G=1$ 
we recover the usual definition
of a tensor category. 

We say that a crossed $G$\nobreakdash-\hspace{0pt}category $\mathcal{C}$ 
is \textit{abelian} if 
\begin{itemize}
\item each component $\mathcal{C}_{\alpha}$ is an abelian category,
\item for all objects $U$ in $\mathcal{C}$ both the functor 
      $U\otimes\_$ and $\_\otimes U$ are both additive and exact, and
\item for all $\beta\in G$, the functor $\varphi_{\beta}$ is
      both additive and exact.
\end{itemize}
Let $\mathcal{D}$ be an abelian tensor category. 
A functor $F\colon\mathcal{C}\to\mathcal{D}$ is \textit{exact} if
each functor $F_{\alpha}\colon\mathcal{C}_{\alpha}\to\mathcal{D}$
is exact.

\begin{rmk}
 Let us recall the definition of the category of crossed 
 $G$\nobreakdash-\hspace{0pt}sets~\cite{FY1,FY2,JS} 
 (see also~\cite{DS}, Example~9), also called
 \textit{$G$\nobreakdash-\hspace{0pt}automorphic sets}~\cite{Brie}
 or \textit{$G$\nobreakdash-\hspace{0pt}racks}~\cite{FR}.
 A \textit{$G$\nobreakdash-\hspace{0pt}set} 
 is a set X together with a left action
 $G\times X\colon(\alpha,x)\mapsto\alpha x$ such that 
 $(\alpha\beta)x=\alpha(\beta x)$ and $1x=x$ for all 
 $\alpha,\beta\in G$, and $x\in X$. A 
 \textit{crossed $G$\nobreakdash-\hspace{0pt}set}
 is a $G$\nobreakdash-\hspace{0pt}set $X$ together with a 
 function $\vert\cdot\vert\colon X\to G$ such that 
 $\vert\alpha x\vert=\alpha\vert x\vert\alpha^{-1}$,
 for all $\alpha\in G$ and $x\in X$. A morphism of crossed 
 $G$\nobreakdash-\hspace{0pt}sets
 $f\colon X\to Y$ is a function such that $f(\alpha x)=\alpha f(x)$ and
 $\bigl\vert f(x)\bigr\vert=\vert x\vert$, for all $\alpha\in G$ and $x\in
 X$. The category $\mathcal{E}$ of crossed
 $G$\nobreakdash-\hspace{0pt}sets 
 is a tensor category with the Cartesian product of $G$-sets
 $X\otimes Y$ and by setting
 $\bigl\vert(x,y)\bigr\vert=\vert x\vert\,\vert y\vert$, for all $x\in X$
 and $y\in Y$. The category $\mathcal{E}$ is braided with
 $c_{X,Y}(x,y)=\bigl(\vert x\vert y, x\bigr)$, for all $x\in X$
 and $y\in Y$. Finally $\mathcal{E}$ is balanced with
 $\theta_{X}(x)=\vert x\vert x$ for all $x\in X$. 
 Crossed $G$\nobreakdash-\hspace{0pt}categories are nothing but 
 monoidal objects in the
 $2$-category $\mathrm{Cat}(\mathcal{E})$ of categories in $\mathcal{E}$
 (the $2$-category of objects in $\mathcal{E}$ introduced by 
 Ehresmann~\cite{E1,E2} is discussed, for instance, 
 in~\cite{Borceux-1}). In particular, the definitions
 of braided (balanced etc.) crossed
 $G$\nobreakdash-\hspace{0pt}category
 provided below agree with the general definitions 
 provided in~\cite{JS} and~\cite{DS}.
\end{rmk}

\subsection{Duals}
A \textit{left autonomous crossed $G$\nobreakdash-\hspace{0pt}category}
$\mathcal{C}=(\mathcal{C},(\cdot)^{\ast})$ is a 
crossed $G$\nobreakdash-\hspace{0pt}category $\mathcal{C}$ 
endowed with a choice of left duals $(\cdot)^{\ast}$
which are compatible with the conjugation in the sense that
if $U^{\ast}$ is a left dual of $U$ via 
$d_{U}\colon U^{\ast}\otimes U\to\tunit$ and
$b_{U}\colon U\otimes U^{\ast}\to\tunit$, then we have
\begin{equation*}
 \Phi_{\beta}(b_{U}) = b_{\Phi_{\beta}(U)}\qquad
 \text{and}\qquad \Phi_{\beta}(d_{U}) = d_{\Phi_{\beta}(U)}
\end{equation*} 
for any $\beta\in G$ and $U\in\mathcal{C}$.
Notice that a duality $(b_{U},d_{U})$ for an object $U$
in $\mathcal{C}_{\alpha}$
induces a duality on all $\Phi_{\alpha}(U)$ compatible 
with the previous equation id and only if the equation
\begin{equation*}
(b_{U},d_{U})=\bigl(\Phi_{\beta}(b_{U}),\Phi_{\beta}(d_{U})\bigr)\label{e:good}
\end{equation*}
is satisfied for all $\beta$ which commute with $\alpha$.

Similarly, it possible to introduce the notion of a 
right autonomous crossed $G$\nobreakdash-\hspace{0pt}category. 
An \textit{autonomous crossed $G$\nobreakdash-\hspace{0pt}category}
is a crossed category that is both left and right autonomous.
When $G=1$ we recover the usual definition of (left/right)
autonomous tensor category~\cite{JS0}.

\subsection{Braiding}
A \textit{braiding} in a crossed $G$\nobreakdash-\hspace{0pt}category 
$\mathcal{C}$ is a family of isomorphisms
\begin{equation*}
  c=\Biggl\{c_{U, V}\in\mathcal{C}\biggl(U\otimes V,
  \Bigl(\lidx{U}{V}\Bigr)\otimes U\biggr)\Biggr\}_{U,V\in\mathcal{C}}
\end{equation*}
satisfying the following conditions.
\begin{subequations}\label{e:braiding}
\begin{itemize}
 \item For all arrows $f\in\mathcal{C}_{\alpha}(U,U')$
       $g\in\mathcal{C}(V,V')$ we have
       \begin{equation}\label{e:braiding-a}
        \Bigl(\bigl(\lidx{\alpha}{g}\bigr)\otimes f\Bigr) \circ c_{U,V} 
        = c_{U',V'} \circ (f \otimes g).
      \end{equation}
 \item For any $U,V,W\in\mathcal{C}$ we have
       \begin{equation}\label{e:braiding-b}
        c_{U\otimes V, W} = a_{\lidx{U\otimes V}{W},U,V} \circ 
        (c_{U,\lidx{V}{W}}\otimes V) \circ 
        a^{-1}_{U,\lidx{V}{W},V} \circ (U\otimes c_{V,W})
       \end{equation}
       \noindent and
       \begin{equation}\label{e:braiding-c}
       c_{U, V\otimes W} = a^{-1}_{\lidx{U\otimes V}{W},U,V} \circ 
                 \Bigl(\bigl(\lidx{U}{V}\bigr)\otimes c_{U,W}\Bigr) \circ 
                 a_{\lidx{U}{V},U,W} \circ (c_{U,V}\otimes W) \circ a^{-1}_{U,V,W}\text{.}
       \end{equation}
   \item For any $U,V\in\mathcal{C}$ and $\beta\in G$ we have
         \begin{equation}\label{e:braiding-d}
            \Phi_{\beta}(c_{U,V}) = c_{\Phi_{\beta}(U),\Phi_{\beta}(V)}\text{.}
         \end{equation}
  \end{itemize}
 \end{subequations}
 
\noindent A crossed $G$\nobreakdash-\hspace{0pt}category endowed 
with a braiding is called a 
\textit{braided crossed $G$\nobreakdash-\hspace{0pt}category.} 
When $G=1$, we recover the usual definition of a braided tensor 
category~\cite{JS}.

\subsection{Twist}
A \textit{twist} in a braided crossed
$G$\nobreakdash-\hspace{0pt}category 
$\mathcal{C}$ is a family of isomorphisms
\begin{equation*}
  \theta=\Bigl\{\theta_{U}\colon U\to\lidx{U}{U}\Bigr\}_{U\in\mathcal{C}}
\end{equation*}
satisfying the following conditions.
\begin{subequations}\label{e:twist}
\begin{itemize}
 \item $\theta$ is \textit{natural}, i.e. for any $f\in\mathcal{C}_{\alpha}(U,V)$ 
       we have
       \begin{equation}\label{e:twist-natural}
        \theta_V \circ f = \bigl(\lidx{\alpha}{f}\bigr) \circ \theta_{U}\text{.}
        \end{equation}
 \item For any $U\in\mathcal{C}_{\alpha}$ and $V\in\mathcal{C}_{\beta}$ 
       we have
       \begin{equation}\label{e:twist-main}
        \theta_{U\otimes V}=c_{\lidx{U\otimes V}{V}, \lidx{U}{U}} 
        \circ c_{\lidx{U}{U},\lidx{V}{V}} \circ 
        (\theta_{U}\otimes\theta_{V})\text{.}
       \end{equation}
  \item For any $U\in\mathcal{C}$ and $\alpha\in G$ we have
        \begin{equation}\label{e:twist-ultra}
         \Phi_{\alpha}(\theta_{U}) = \theta_{\Phi_{\alpha}(U)}\text{.}
        \end{equation}
\end{itemize}
\end{subequations}
A braided crossed $G$\nobreakdash-\hspace{0pt}category endowed with a twist
is called a \textit{balanced crossed $G$\nobreakdash-\hspace{0pt}category.}
When $G=1$ we recover the usual definition of 
a balanced tensor category~\cite{JS}.

A \textit{ribbon crossed $G$\nobreakdash-\hspace{0pt}category} $\mathcal{C}$ is a 
balanced crossed $G$\nobreakdash-\hspace{0pt}category that is also 
a left autonomous crossed $G$\nobreakdash-\hspace{0pt}category such that for any 
$U\in\mathcal{C}_{\alpha}$ (with $\alpha\in G\/$), 
\begin{equation}\label{e:tortility} 
\Bigl(\bigl(\lidx{U}{U}\bigr)\otimes\theta_{\leftidx{^{U}}{\! U}{^{\ast}}}\Bigr) \circ 
 b_{\lidx{U}{U}}=
 (\theta_{U} \otimes U^{\ast}) \circ b_{U}\text{.}
\end{equation}
When $G=1$ we recover the usual definition of a ribbon 
category~\cite{RT,Tur-QG}, also called 
\textit{tortile tensor category}~\cite{JS,JS-tortile,Shum}.

\section{Crossed Hopf Group Algebras}
 We recall the definition of a crossed Hopf group algebra as
 in~\cite{Zunino-2} and we introduce the notion of a cobraided 
 and a coribbon crossed Hopf group algebra.\medskip
 
\subsection{Basic Definitions}
 A crossed Hopf $G$\nobreakdash-\hspace{0pt}algebra $H$ is a family 
 $\bigl\{(H_{\alpha},\Delta_{\alpha},\eta_{\alpha})\bigr\}_{\alpha\in G}$
  of coalgebras 
 endowed with the following data.
 \begin{itemize}
  \item A family of coalgebra morphisms
        $\mu_{\alpha,\beta}\colon H_{\alpha}\otimes
        H_{\beta}\to H_{\alpha\beta}$,
        the \textit{multiplication,} that is 
        associative in the sense that we have
        \begin{equation}\label{e:associativity}
         \mu_{\alpha\beta,\gamma}\circ (\mu_{\alpha,\beta}\otimes H_{\gamma})=
         \mu_{\alpha,\beta\gamma}\circ (H_{\alpha}\otimes\mu_{\beta,\gamma})\text{.}
        \end{equation}
        If, for all $h\in H_{\alpha}$ and $k\in H_{\beta}$, we set
        $hk=\mu_{\alpha,\beta}(h,k)$, then~\eqref{e:associativity} 
        can be simply rewritten as the usual associativity law $(hk)l = h(kl)$
        for all $h\in H_{\alpha}$, $k\in H_{\beta}$, $l\in H_{\gamma}$.
  \item An algebra morphism $\eta\colon \Bbbk\to  H_{1}$,
        the \textit{unit,} such that, if we set 
        $1=\eta(1_{\Bbbk})$, then $1h = h = h1$ for any $h\in
        H_{\alpha}$ and $\alpha\in G$.
  \item A set of coalgebra isomorphisms 
        $\varphi_{\beta}=\colon H_{\alpha}\to H_{\beta\alpha\beta^ {-1}}$,
        the \textit{conjugation,} such that 
        \begin{itemize}
         \item $\varphi_{\beta}(hk) = \varphi_{\beta}(h)\varphi_{\beta}(k)$,
         \item $\varphi_{\beta}(1) = 1$.
        \end{itemize}
  \item A set of linear isomorphisms  
        $S_{\alpha} \colon H_{\alpha}\to H_{\alpha^{-1}}$,
        the \textit{antipode,} such that
        \begin{equation*}
        \mu_{\alpha^{-1},\alpha}\circ (S_{\alpha}\otimes H_{\alpha})\circ \Delta_{\alpha}=
        \eta\circ\varepsilon_{\alpha}=\mu_{\alpha,\alpha^{-1}}\circ 
        (H_{\alpha}\otimes S_{\alpha})\circ \Delta_{\alpha}\text{.}
        \end{equation*}
 \end{itemize}
\noindent The coalgebra $H_{\alpha}$ is called the \textit{\alphath component of $H$.}
We say that $H$ is \textit{of finite type} if $\dim H_{\alpha}<\infty$ for all
$\alpha\in G$.

\subsection{Packed Form of a Crossed Group Algebra}\label{par:packed-alg}
 An equivalent definition of crossed Hopf $G$\nobreakdash-\hspace{0pt}algebra can be obtained as
 follows.
 
 Let $H$ be a crossed Hopf $G$\nobreakdash-\hspace{0pt}algebra. We obtain a Hopf algebra 
 $H_{\pk}$, which we call the \textit{packed form of $H$.}
 As a $G$\nobreakdash-\hspace{0pt}graded coalgebra, $H_{\pk}$ is the direct sum of the components of $H$. 
 Let $i_{\alpha}$ be the inclusion of $H_{\alpha}\hookrightarrow H_{\pk}$.
 The multiplication $\mu_{\pk}$ of $H_{\pk}$ is, by definition, the
 colimit of linear morphisms 
 $\mbox{$\varinjlim$}_{\alpha,\beta\in G}(i_{\alpha\beta}\circ\mu_{\alpha,\beta})$, 
 i.e., the unique linear map from $H_{\pk}\otimes H_{\pk}$ to $H_{\pk}$ 
 such that the restriction to $H_{\alpha}\otimes H_{\beta}\subset H_{\pk}\otimes H_{\pk}$ 
 coincides with $\mu_{\alpha,\beta}$. 
 The unit $H_{\pk}$ is $1_{\pk}=1\in H_{1}\subset H_{\pk}$. The antipode of $H_{\pk}$ 
 is the sum $S_{\pk}=\sum_{\alpha\in G}S_{\alpha}$. The Hopf algebra $H_{\pk}$ is 
 endowed with a group morphism 
 $\map{\varphi}{G}{\aut(H_{\pk})}{\alpha}{\varphi_{\pk,\alpha}}$, 
 where $\varphi_{\pk,\alpha}=\sum_{\beta\in G}
 \varphi_{\alpha}^{\beta}\colon\bigoplus_{\beta\in G}H_{\beta}
 \longmapsto\bigoplus_{\beta\in G}H_{\beta}$.
 
 Conversely, let $H_{\grande}=(H_{\grande},\varphi_{\grande})$ 
 be a Hopf algebra (with 
 multiplication $\mu_{\grande}$, unit $1$, and antipode $S_{\grande}\/$)
 endowed with a group homomorphism
 $\map{\varphi_{\grande}}{G}{\aut(H_{\grande})}{\alpha}{\varphi_{\grande,\alpha}}$.
 Suppose that the following conditions are satisfied.
 \begin{itemize}
  \item There exists a family of sub-coalgebras
        $\{H_{\alpha}\}_{\alpha\in G}$ of  $H_{\grande}$
        such that $H_{\grande}=\bigoplus_{\alpha\in G}H_{\alpha}$.
  \item $H_{\alpha}\cdot H_{\beta}\subset H_{\alpha\beta}$ for all
        $\alpha,\beta\in G$.
  \item $1\in H_{1}$.
  \item $\varphi_{\grande,\beta}$ sends $H_{\alpha}\subset H_{\grande}$ to
        $H_{\beta\alpha\beta^{-1}}\subset H_{\grande}$.
  \item $S_{\grande}(H_{\alpha})\subset H_{\alpha^{-1}}$.
 \end{itemize}
\noindent In the obvious way $H_{\grande}$ gives rise to a 
crossed Hopf $G$\nobreakdash-\hspace{0pt}algebra $H$ such that 
$H_{\pk}=H_{\grande}$.
We say that $H_{\pk}$ is the \textit{packed form of $H$.}

\begin{rmk}
By dualizing in obvious way the axioms for a crossed Hopf
$G$\nobreakdash-\hspace{0pt}algebra
one get the dual notion of crossed Hopf
$G$\nobreakdash-\hspace{0pt}coalgebra, see~\cite{Tur-CPC}.
In particular, if $H$ is a crossed Hopf
$G$\nobreakdash-\hspace{0pt}algebra of finite type
(i.e. $\dim H_{\alpha}<+\infty$ for all $\alpha\in G\/$),
then, by taking the linear duals of the components of $H$,
we get in obvious way a crossed Hopf
$G$\nobreakdash-\hspace{0pt}coalgebra $H^{\ast}$,
the \textit{dual of $H$.}
Notice, however, that, in general, it is not possible 
to describe a crossed Hopf 
$G$\nobreakdash-\hspace{0pt}coalgebra
as a graded Hopf algebra.
\end{rmk}

We say that a crossed Hopf $G$\nobreakdash-\hspace{0pt}algebra $H$
is \textit{cosemisimple} if every component of $H$ is cosemisimple 
as a coalgebra. It is proved in~\cite{Virelizi} that, if $K$ is a
finite type crossed Hopf $G$\nobreakdash-\hspace{0pt}algebra, then
the components of $K$ are semisimple algebras in and only if the
component $K_{1}$ is a semisimple Hopf algebra. Therefore, by duality,
we obtain that a crossed Hopf $G$\nobreakdash-\hspace{0pt}algebra 
$H$ of finite type
is cosemisimple if and only if its component $H_{1}$ is a cosemisimple
Hopf algebra.

\subsection{Categories of Comodules}
Let $H$ be a crossed Hopf $G$\nobreakdash-\hspace{0pt}algebra.
We define the 
crossed $G$\nobreakdash-\hspace{0pt}category $\Com H$ of left
$H$\nobreakdash-\hspace{0pt}comodules as follow.
\begin{itemize}
 \item $\Com_{\alpha}H$ is the category of 
       $H_{\alpha}$\nobreakdash-\hspace{0pt}comodules.
 \item The tensor unit $\tunit$ is the ground field $\C$
       with the comodule structure given by the the unit
       of $H$ (i.e. $c\mapsto c\otimes 1_{H}$ for all scalars
       $c\/$).
 \item Let $M$ be a $H_{\alpha}$\nobreakdash-\hspace{0pt}comodule
       with coaction
       $\map{\Delta_{M}}{M}{H_{\alpha}\otimes M}%
        {m}{m_{\alpha}\otimes m_{M}}$
       and let $N$ be a $H_{\beta}$\nobreakdash-\hspace{0pt}comodule
       with coaction
       $\map{\Delta_{N}}{N}{H_{\beta}\otimes N}%
        {n}{n_{\beta}\otimes n_{N}}$. Define $M\otimes N$ as the vector
       space $M\otimes_{\C}N$ with the  
       $H_{\alpha\beta}$\nobreakdash-\hspace{0pt}comodule structure
       given by 
       \begin{equation*}
       \Delta_{M\otimes N}\colon
       m\otimes n\mapsto m_{\alpha}n_{\beta}\otimes m_{M}\otimes n_{N}\text{.}
       \end{equation*}
 \item The functors $\Phi_{\beta}$ are defined as follows.
       Let $M$ be a $H_{\alpha}$\nobreakdash-\hspace{0pt}comodule.
       The 
       $H_{\beta\alpha\beta^{-1}}$\nobreakdash-\hspace{0pt}comodule
       $\lidx{\beta}{M}=\Phi_{\beta}(M)$ is isomorphic to $M$ as a vector space
       and the element in $\lidx{\beta}{M}$ corresponding to
       $m\in M$ is denoted $\lidx{\beta}{m}$. The coaction of 
       $H_{\beta\alpha\beta^{-1}}$ is obtained by setting
       \begin{equation*}
       (\lidx{\beta}{m})_{\beta\alpha\beta^{-1}}\otimes
       (\lidx{\beta}{m})_{\lidx{\beta}{M}}=
       \varphi_{\beta}(m_{\alpha})\otimes\lidx{\beta}{(m_{M})}\text{.}
       \end{equation*}
       
       If $f\colon M\to N$ is a morphism of
       $H_{\alpha}$\nobreakdash-\hspace{0pt}comodules,
       then set
       \begin{equation*}\tmap
        {\lidx{\beta}{f}=\Phi_{\beta}(f)}
        {\lidx{\beta}(M)}{\lidx{\beta}(N)}
        {\lidx{\beta}{m}}{\lidx{\beta}{\bigl(f(m)\bigr)}}\text{.}
       \end{equation*}
\end{itemize}

\begin{lemma}\label{l:indent}
For all $\alpha,\beta\in G$ the map
\begin{equation*}
 \map{\hat{\varphi}_{\beta}}
 {\lidx{\beta}{(H_{\alpha})}}
 {H_{\beta\alpha\beta^{-1}}}
 {\lidx{\beta}{h}}{\varphi_{\beta}(h)}
\end{equation*}
is an isomorphism of $H_{\beta\alpha\beta^{-1}}$\nobreakdash-\hspace{0pt}comodules.
\end{lemma}

The proof is straightforward and left to the reader.

 A \textit{cobraided crossed Hopf $G$\nobreakdash-\hspace{0pt}algebra}
 is a crossed Hopf $G$\nobreakdash-\hspace{0pt}algebra $H$ endowed
 with a family of linear maps
 \begin{equation*}
  \gamma_{H}=\{\gamma_{\alpha,\beta}\colon H_{\alpha}\otimes
  H_{\beta}\to \C\}
 \end{equation*}
 satisfying the following conditions.
 \begin{itemize}
 \begin{subequations}\label{e:gamma}
  \item There exists a family of linear maps
        $\tilde{\gamma}_{\alpha,\beta}\colon H_{\alpha}\otimes H_{\beta}\to\C$
        such that, for all $h\in H_{\alpha}$ and $k\in H_{\beta}$,
        \begin{equation}
         \tilde{\gamma}_{\alpha,\beta}(h'\otimes k')
         \gamma_{\alpha,\beta}(h''\otimes k'')=
         \gamma_{\alpha,\beta}(h'\otimes k')
         \tilde{\gamma}_{\alpha,\beta}(h'\otimes k')=
         \varepsilon(h)\varepsilon(k)\text{.}
        \end{equation}
  \item For all $h\in H_{\alpha}$ and $k\in H_{\beta}$, we have
        \begin{equation}
         k'h'\gamma_{\alpha,\beta}(h''\otimes k'')=
         \gamma_{\alpha,\beta}(h'\otimes k')
         \varphi_{\beta}(h'')k''\text{.}
        \end{equation}
  \item For all $h_{1}\in H_{\alpha_{1}}$, $h_{2}\in H_{\alpha_{2}}$, and
        $k\in H_{\beta}$, we have
        \begin{equation}
         \gamma_{\alpha_{1}\alpha_{2},\beta}\bigl((h_{1}h_{2})\otimes
         k\bigr)=
         \gamma_{\alpha_{1},\beta}(h_{1}\otimes k')
         \gamma_{\alpha_{2},\beta}(h_{2}\otimes k'')\text{.}
        \end{equation}
  \item For all $h\in H_{\alpha}$, $k_{1}\in H_{\beta_{1}}$, and
        $k_{2}\in H_{\beta_{2}}$, we have
        \begin{equation}
         \gamma_{\alpha,\beta_{1}\beta_{2}}\bigl(h\otimes
         (k_{1}k_{2})\bigr)= 
         \gamma_{\beta_{2}\alpha\beta_{2}^{-1},\beta_{1}}
         \bigl(\varphi_{\beta_{2}}(h'')\otimes k_{1}\bigr)
         \gamma_{\alpha,\beta_{2}}(h'\otimes k_{2})\text{.}
        \end{equation}
  \item For all $h\in H_{\alpha}$, $k\in H_{\beta}$, and $\lambda\in G$, we have
        \begin{equation}
         \gamma_{\lambda\alpha\lambda^{-1},\lambda\beta\lambda^{-1}}
         \bigl(\varphi_{\lambda}(h)\otimes\varphi_{\lambda}(k)\bigr)=
         \gamma_{\alpha,\beta}(h\otimes k)\text{.}
        \end{equation}
 \end{subequations}
 \end{itemize}

\begin{thm}\label{thm:Beethoven}
There is a \textup{1--1} correspondence between braidings in $\Com H$ 
\textup{(}as a crossed
$G$\nobreakdash-\hspace{0pt}category\/\textup{)}
and cobraided crossed Hopf $G$\nobreakdash-\hspace{0pt}algebras
structures on $H$.
\end{thm}

\begin{proof}
 Let $H$ be a cobraided crossed Hopf $G$\nobreakdash-\hspace{0pt}algebra.
 We obtain a braiding in $\Com H$ by setting, 
 for every $H_{\alpha}$-\hspace{0pt}comodule $M$ and every 
 $H_{\beta}$-\hspace{0pt}comodule $M$,
 \begin{equation*}
  c^{\gamma}_{M,N}(m\otimes n)=\gamma_{\beta,\alpha}(n_{\beta}\otimes
  m_{\alpha})\Bigl(\lidx{\alpha}{(n_{N})}\Big)\otimes m_{M}\text{.}
 \end{equation*}
 Conversely, let $c$ be a braiding in $\Com H$. 
 We obtain a cobraided structure of $H$ by setting, for all $\alpha,\beta\in G$,
 \begin{equation*}
  \gamma^{c}_{\alpha,\beta}\colon H_{\alpha}\otimes H_{\beta}
  \xrightarrow{c_{H_{\alpha,H_{\beta}}}}
  \Bigl(\lidx{\alpha}{H_{\beta}}\Bigr)\otimes H_{\alpha}
  \xrightarrow{\hat{\varphi}_{\alpha}\otimes H_{\alpha}}
  H_{\alpha\beta\alpha^{-1}}\otimes H_{\alpha}
  \xrightarrow{\varepsilon_{\alpha\beta\alpha^{-1}}\otimes\varepsilon_{\alpha}}
  \C\text{,}
 \end{equation*}
 where $\hat{\gamma}_{\alpha}$ is defined as in Lemma~\ref{l:indent}.
 The proof of the theorem is an adaptation of the standard proof for
 Hopf algebras (see~\cite{Kas} and~\cite{JS-Tannaka}).
\end{proof}

A \textit{coribbon crossed Hopf $G$\nobreakdash-\hspace{0pt}algebra}
is a crossed Hopf $G$\nobreakdash-\hspace{0pt}algebra $H$ endowed with a
cobraided strcuture $\gamma$ and a family
\begin{equation*}
 \tau_{H}=\{\tau_{\alpha}\colon H_{\alpha}\to\C\}\text{,}
\end{equation*}
called \textit{cotwisted structure,} satisfying the following conditions.
\begin{subequations}\label{e:tau}
\begin{itemize}
\item There exists a family of linear maps 
      $\tilde{\tau}_{\alpha}\colon H_{\alpha}\to\C$ such that, for all
      $h\in H_{\alpha}$,
      \begin{equation*}
       \tau_{\alpha}(h')\tilde{\tau}_{\alpha}(h'')=
       \tilde{\tau}_{\alpha}(h')\tau_{\alpha}(h'')=
       \varepsilon_{\alpha}(h)\text{.}
      \end{equation*}
\item For all $h\in H_{\alpha}$ we have
      \begin{equation*}
       \tau_{\alpha}(h')\varphi_{\alpha}(h'')=h'\tau_{\alpha}(h'')\text{.}
      \end{equation*}
\item For all $h\in H_{\alpha}$ and $k\in H_{\beta}$, we have
      \begin{equation*}
       \tau_{\alpha\beta}(hk)=
       \gamma_{\beta,\alpha}(k'_{\beta}\otimes h'_{\alpha})
       \tau_{\alpha}(h''_{\alpha})
       (\tau_{\alpha\beta\alpha^{-1}}\circ\varphi_{\alpha})(k''_{\beta})
       (\varphi_{\alpha}\circ\gamma_{\alpha,\beta})
       \bigl(h'''_{\alpha}\otimes\varphi_{\beta}(k'''_{\beta})\bigr)\text{.}
      \end{equation*}
\item For all $\alpha\in G$, we have
      \begin{equation*}
       \tau_{\alpha^{-1}}\circ S_{\alpha}=\tau_{\alpha}\text{.}
      \end{equation*}
\item For all $\alpha,\beta\in G$,
      \begin{equation*}
       \tau_{\beta\alpha\beta^{-1}}\circ\varphi_{\beta}=\tau_{\alpha}\text{.}
      \end{equation*}
\end{itemize}
\end{subequations}

\begin{thm}
There is a \textup{1--1} correspondence between twist in $\Com H$ and cotwisted
structures in $H$.
\end{thm}

\begin{proof}
Let $H$ be a cotwisted crossed Hopf $G$\nobreakdash-\hspace{0pt}algebra.
We obtain a twist in $\Com H$ by setting, 
for every $H_{\alpha}$\nobreakdash-\hspace{0pt}comodule $M$,
\begin{equation*}
 \theta^{\tau}_{M}(m)=\tau_{\alpha}(m_{\alpha})\lidx{\alpha}{(m_{M})}\text{.}
\end{equation*}
Conversely, let $\theta$ be a twit in $\Com H$.  
We obtain a cotwisted structure of $H$ by setting, for all $\alpha\in G$,
\begin{equation*}
\tau^{\theta}_{\alpha}\colon
H_{\alpha}\xrightarrow{\theta_{H_{\alpha}}}\leftidx{^{\alpha}}{H}{_{\alpha}}
\xrightarrow{\hat{\varphi}_{\alpha}}H_{\alpha}
\xrightarrow{\varepsilon_{\alpha}}\C\text{.}
\end{equation*}
Again, the proof of the theorem is an adaptation of the standard one.
\end{proof}

\subsection{Modularity}
We give a definition of modular coribbon Hopf algebra
by transposing the usual axioms of a modular
(ribbon) Hopf algebra, see e.g.~\cite{Tur-QG} or~\cite{Kas}. 
Then we introduce the notion of modular coribbon Hopf 
$G$\nobreakdash-\hspace{0pt}algebra.

Let $H_{1}$ be a coribbon Hopf algebra, let $U$ be a 
finite-dimensional $H_{1}$\nobreakdash-\hspace{0pt}comodule, 
and let $f\colon U\to U$ an endomorphism of
$H_{1}$\nobreakdash-\hspace{0pt}comodules. We define the
\textit{quantum trace of f} as the scalar
\begin{equation*}
\qtr(f)=
\sum_{i=1}^{n}\tau_{1}\bigl((e_{i})_{H_{1}}\bigr)
\tilde{\gamma}_{1,1}\bigl((e^{i})_{H_{1}}\otimes(e_{i})''_{H_{1}}\bigr)
f\bigl((e_{i})_{U}\bigr)(e^{i})_{U^{\ast}}\text{,}
\end{equation*}
being $\{e_{i}\}_{i=1}^{n}$ a basis of $U$, $\{e^{i}\}_{i=1}^{n}$
the dual basis
$\tilde{\gamma}_{1,1}$ the inverse of the cobraiding
of $H_{1}$, and $\tau_{1}$ the coribbon of $H_{1}$
(notice that this is nothing but the usual quantum
trace of $f$ inside the category $\Com H_{1}\/$). 
We say that $U$ is \textit{negligible} if $\qtr(\Id_{U})=0$.

We say that a coribbon Hopf algebra $H_{1}$ is \textit{modular}
if it is endowed with a finite family of simple finite-dimensional
$H_{1}$\nobreakdash-\hspace{0pt}comodules $\{V_{i}\}_{i\in I}$ 
satisfying the following conditions.
\begin{itemize}
\item There exists an element $0\in I$ such that $V_{0}=\C$
      (with the structure of $H_{1}$\nobreakdash-\hspace{0pt}comodule
      given by the multiplication).
\item For any $i\in I$, there exists $i^{\ast}\in I$ such that
      $V_{i^{\ast}}$ is isomorphic to $V^{\ast}_{i}$.
\item For any $j,k\in I$, the $H_{1}$\nobreakdash-\hspace{0pt}comodule
      $V_{j}\otimes V_{kj}$ is isomorphic to a finite sum of elements
      of $\{V_{i}\}_{i\in I}$, possibly with repetitions, and a
      negligible $H_{1}$\nobreakdash-\hspace{0pt}comodule.
\item Let $\mathfrak{s}_{ij}$ be the quantum trace of the 
      endomorphism
      \begin{equation*}
       \tmap{c_{V_{j^{\ast},i},V_{i}}\circ c_{V_{i},V_{j^{\ast}}}}
           {V_{i}\otimes V_{j^{\ast}}}{V_{i}\otimes V_{J^{\ast}}}
           {x\otimes y}%
           {\gamma_{1,1}(y'_{H_{1}}\otimes x'_{H_{1}})
            \gamma_{1,1}(x''_{H_{1}}\otimes y''_{H_{1}})
            x_{V_{i}}\otimes y_{V_{j^{\ast}}}}\text{.}
      \end{equation*}
      The matrix $(\mathfrak{s}_{i,j})_{i,j\in I}$ is invertible.
\end{itemize}
A coribbon Hopf $G$\nobreakdash-\hspace{0pt}algebra $H$ is \textit{modular}
if its component $H_{1}$ is modular. We observe that, 
if $\dim H_{1}<+\infty$, then $H_{1}$ is modular if and only if
$H_{1}^{\ast}$ is modular in the usual sense. Therefore, if $H$ is
of finite type, then it is modular if and only if $H^{\ast}$ is modular
in the usual sense~\cite{Tur-CPC,Virelizi}. Also, $H_{1}$ is modular
if and only if the category $\Com H_{1}$ is modular in the usual sense.

\section{Tannaka Theory for Tensor Categories}
We recall some basic facts about Tannaka Theory. For a classical
reference, see~\cite{SR}. For a reference about Tannaka Theory
and braided tensor categories see~\cite{JS-Tannaka}, of which
we follows the approach, or~\cite{Sch}. In a couple of cases we also
reproduce a short sketch of the proof since the notation introduced
will be used in the sequel.\medskip

\subsection{Dinatural Transformations} Let both $\mathcal{C}$
and $\mathcal{B}$ be two small categories and let both $S$  and $T$ be 
functors from $\Op{\mathcal{C}}\times\mathcal{C}$ to $\mathcal{B}$.
A \textit{dinatural transformation} $\alpha\colon
S\xrightarrow{\cdot\cdot}T$ (see,~e.g.,~\cite{McL}) 
is a function that assigns to an object $C\in\mathcal{C}$
an arrow $\alpha_{C}\colon S(C,C)\to T(C,C)$ in $\mathcal{B}$
such that, for all arrows $f\colon C\to D$ in $\mathcal{C}$,
the diagram
  \begin{equation}\label{e:T1}
     \vcenter{\xymatrix{\ &  S(C,C)\ar[r]^{\alpha_{C}} & 
     T(C,C)\ar@/^1pc/[dr]^(.65){T(C,f)} & \ \\
     S(D,C)\ar@/^1pc/[ur]^(.35){S(f,C)}\ar@/_1pc/[dr]_(.35){S(D,f)} & \ & \  & T(C,D) \\
     \ & S(D,D)\ar[r]_{\alpha_{D}} &
     T(D,D)\ar@/_1pc/[ur]_(.65){T(f,D)} & \ }}\text{.}
   \end{equation}
commutes.

Let $B$ be an object of $\mathcal{B}$, and $S=B$, i.e.\@ the
constant functor fixed by $B$. We can rewrite the commutativity of~\eqref{e:T1} as
  \begin{equation*}
  T(C,f)\circ\alpha_{C}=T(f,D)\circ\alpha_{D}\text{.}
  \end{equation*}
We say that $B$ is an \textit{end of $T$} 
if $\alpha$ is universal for the property that, for any
dinatural transformation 
$\alpha'\colon B'\xrightarrow{\cdot\cdot}T$,
there exists an arrow $b\in\mathcal{B}(B',B)$ such that the diagram
  \begin{equation*}
   \vcenter{\xymatrix @C=2pc{B &\ & \ B'\ar@/_1.5pc/[ll]_{b}
   \ar@/^/[dl]^(.35){\alpha'}\\ & T(C,C)\ar@{<-}@/^/[ul]^(.65){\alpha} &  \ }}
  \end{equation*}
commutes for all objects $C\in\mathcal{C}$. If $T$ has an end, then
this end is unique up to canonical isomorphism.

Let us consider again diagram~\eqref{e:T1} and suppose this time 
that $T=B$ for a fixed object $B\in\mathcal{C}$.  We can rewrite the
commutativity of~\eqref{e:T1} as
  \begin{equation*}
   \alpha_{C}\circ S(f,C)=\alpha_{D}\circ S(D,f)\text{.}
  \end{equation*}
We say that $B$ is a \textit{coend of $S$} if $\alpha$ 
is universal for the property that for any $\alpha'\colon S\xrightarrow{\cdot\cdot}B'$,
there exists an arrow $b\in\mathcal{B}(B, B')$ such that the diagram 
  \begin{equation*}
   \vcenter{\xymatrix @C=2pc{B &\ & \ B'\ar@{<-}@/_1.5pc/[ll]_{b}
   \ar@{<-}@/^/[dl]^(.35){\alpha'}\\ & S(C,C)\ar@/^/[ul]^(.65){\alpha} &  \ }}
  \end{equation*}
commutes for all $C\in\mathcal{C}$. Again, if $S$ has a
coend, then this coend is unique up to canonical isomorphism.

\subsection{Ends and Coends as Vector Spaces}
Let both $X$ and $Y$ be functors from a small category $\mathcal{C}$ to
the category $\Vect$ of complex vector spaces.

\begin{lemma}
The bifunctor $\Hom_{\Bbbk}\bigl(X(\dummy[1]),Y(\dummy[2])\bigr)$ 
has an end $\Hom(X,Y)$.
\end{lemma}

\begin{sketch}
To construct $\Hom(X,Y)$, for every morphism $f\colon A\to B$
in $\mathcal{C}$ define two linear maps
 \begin{equation*}
   p_{f}, q_{f}\colon
   \prod_{C\in\mathcal{C}}\Hom_{\C}\bigl(X(C),Y(C)\bigr)
   \longmapsto\Hom_{\Bbbk}\bigl(X(A),Y(B)\bigr)
  \end{equation*}
by
  \begin{equation*}\begin{split}
      p_{f}\colon u=\{u_{C}\vert C\in\mathcal{C}\}\mapsto Y(f)\circ u_{A}\text{,}\\
      q_{f}\colon u=\{u_{C}\vert C\in\mathcal{C}\}\mapsto u_{B}\circ X(f)\text{,}
  \end{split}\end{equation*}
and $\Hom(X,Y)$ as the equalizer of all pairs $(p_{f},q_{f})$.
When $X=Y$ we will use the notation $\End(X)=\Hom(X,Y)$.
\end{sketch}

\begin{lemma}
If both $X$ and $Y$ send all objects to finite dimensional vector
spaces, then the bifunctor
$\Hom_{\Bbbk}\bigl(X(\dummy[1]),Y(\dummy[2])\bigr)$ has a coend 
$\HomCheck(X,Y)$ and $\Hom(X,Y)$ is the dual of $\HomCheck(X,Y)$.
\end{lemma}

\begin{sketch}
Define $\HomCheck(X,Y)$ as the coequalizer of all pairs
$\bigl(\leftidx{^{t}}{p}{_{f}},\leftidx{^{t}}{\! q}{_{f}}\bigr)$,
where  $f\colon A\to B$ is an arrow in $\mathcal{C}$ and
\begin{equation*}
    \leftidx{^{t}}{p}{_{f}},\leftidx{^{t}}{\! q}{_{f}}\colon \Hom_{\C}\bigl(X(A),Y(B)\bigr)^{\ast}\longmapsto\sum_{C\in\mathcal{C}}
    \Hom_{\C}\bigl(X(C),Y(C)\bigr)^{\ast}\text{.}
\end{equation*}
When $X=Y$ we will use the notation $\EndCheck(X)=\Hom(X,Y)$.
\end{sketch}

Let $U$ and $V$ be two finite dimensional vector spaces. By means of the
canonical isomorphisms
$\Hom_{\C}(U,V)^{\ast}\cong V^{\ast}\otimes_{\C} U\cong\Hom_{\C}(V,U)$
we have a pairing
\begin{equation*}
    \map{\langle\dummy[1],\dummy[2]\rangle}{\Hom_{\C}(U,V)\times\Hom_{\C}(V,U)}{\C}{(h,k)}{\Tr(hk)}\text{.}
\end{equation*} 
By using this, $\HomCheck(X,Y)$ can be defined as the common
coequalizer of all maps
\begin{equation*}
 i_{f},\,j_{f}\colon\Hom\bigl(Y(B),X(A)\bigr)\to
 \sum_{C\in\mathcal{C}}\Hom\bigl(Y(C),X(C)\bigr)\text{,}
\end{equation*}
where, for any $h\in\Hom_{\C}\bigl(Y(B),X(A)\bigr)$, we set
  \begin{align*}
      i_{f}(h) &=\bigl(h\circ Y(f), A\bigr)\ \ \text{and}\\
      j_{f}(h) &=\bigl(X(f)\circ h, B\bigr)\text{.}
  \end{align*}
Here the second component of a pair $(h\circ Y(f), A\bigr)$
indicates to which component of the direct sum this element belongs.

For any object $C$ in $\mathcal{C}$ and any map
$h\in\Hom_{\C}\bigl(Y(C),X(C)\bigr)$, let
$[h]$ be the image of $h$ under the canonical map
$\Hom_{\C}\bigl(Y(C),X(C)\bigr)\to\HomCheck(X,Y)$. The space
$\HomCheck(X,Y)$ is generated as a vector space by the symbols $[h]$
subject to the relations
\begin{itemize}
 \item $[c_{1} h_{1} + c_{2} h_{2}]=c_{1}[h_{1}]+ c_{2}[h_{2}]$ for all
       $h_{1},h_{2}\in\Hom\bigl(Y(C),X(C)\bigr)$, and
       $c_{1},c_{2}\in\C$,\smallskip
 \item $\bigl[k\circ Y(f)\bigr]=\bigl[X(f)\circ k\bigr]$ 
       for all $f\colon A\to B$ in $\mathcal{C}$ and 
       $k\in\Hom_{\C}\bigl(Y(B),X(A)\bigr)$.
\end{itemize}
The pairing between $\Hom(X,Y)$ and $\HomCheck(X,Y)$ is given by
\begin{equation*}
    \map{\langle\dummy[1],\dummy[2]\rangle}%
    {\Hom(X,Y)\times\HomCheck(X,Y)}{\C}{\bigl(u,[h]\bigr)}%
    {\Tr(u_{C}\circ h)}\text{,}
\end{equation*}
where $h\in\Hom_{\C}\bigl(Y(C),X(C)\bigr)$.

\subsection{The Coalgebra $\EndCheck(X)$}
Let $X\colon\mathcal{C}\to\Vect$ be a functor whose values are finite
dimensional vector spaces. The space $\EndCheck(X)=\HomCheck(X,X)$
is a coalgebra as follows.

Recall that for any vector space $V$ the space $\End(V)$ is a
coalgebra via
\begin{equation*}
 \map{\delta}{\End(V)}{\End(V)\otimes\End(V)}{e^{i}_{j}}{e^{i}_{k}\otimes e^{k}_{j}}\text{,}
\end{equation*}
where $e_{1}, \ldots, e_{n}$ is a basis of $V$,
$e^{i}_{j}=e_{i}^{\ast}\otimes e_{j}$, and on the right hand side the
sum runs on $k=1, \ldots, n$. The counit of $\End(V)$ is the trace 
$\Tr\colon\End(V)\to\C$.

The coalgebra structure on $\EndCheck(X)$ is given by
\begin{equation*}
 \Delta[e^{i}_{j}]=[e^{i}_{k}]\otimes [e^{k}_{j}]\text{,}
\end{equation*}
with counit
\begin{equation*}
 \varepsilon[h]=\Tr(h)\text{,}
\end{equation*}
that is, $\EndCheck(X)$ is a quotient of direct sums of coalgebras via
the canonical map
\begin{equation*}
 \sum_{C\in\mathcal{C}}\End\bigl(X(C)\bigr)\to\EndCheck(X)\text{.}
\end{equation*}

The following lemma will be crucial in the next section. Let both
$X\colon\mathcal{C}\to\Vect$ and 
$Y\colon\mathcal{D}\to\Vect$ be functors whose values are finite
dimensional vector spaces. Set 
\begin{equation*}
 X\otimes Y\colon \mathcal{C}\times\mathcal{D}\to\Vect\colon
                  (U,V)\to U\otimes V\text{.}
\end{equation*}

\begin{lemma}\label{l:noja}
The map
\begin{equation*}
\dmap{m}{\EndCheck(X)\otimes\EndCheck(Y)}{\EndCheck(X\otimes Y)}%
        {[S]\otimes [T]}{[S\otimes T]}
\end{equation*}
\textup{(\/}for any object $A$ in $\mathcal{C}$,
and $B$ in $\mathcal{D}$, and any $S\in\End\bigl(X(A)\bigr)$
and $T\in\End\bigl(Y(B)\bigr)$\textup{)}
is a coalgebras isomorphism.
\end{lemma}

Let both $X$ and $Y$ as above. We say that 
$F\colon\mathcal{C}\to\mathcal{D}$ is an \textit{equivalence}
[resp.\@ \textit{an isomorphism\/}] of Tannakian categories
it is an equivalence [resp.\@ an isomorphism] of categories
and $X=Y\circ F$.

\begin{lemma}
The map $\EndCheck(X)\to\EndCheck(Y):[h]\to[h]$ 
\textup{(}for all objects $C$ in $\mathcal{C}$ and all 
$h\in\End\bigl(X(C)\bigr)=
\End\bigl((Y\circ F)(C)\bigr)\/$\textup{)}
is a morphism of coalgebras
and it is an isomorphism if and only if
$F$ is an equivalence of Tannakian categories.
\end{lemma}

We also need a small variant of the previous lemma.

\begin{lemma}\label{l:trivial}
Let $F\colon\mathcal{C}\to\mathcal{D}$ be an isomorphism if categories
and let $f\colon X\to Y\circ F$ be a natural isomorphism. The map
$\EndCheck(X)\to\EndCheck(Y)\colon[h]\to[f_{C}\circ h\circ f^{-1}_{C}]$
is an isomorphism of coalgebras.
\end{lemma}

The proof is straightforward and left to the reader.

The main theorem of Tannaka Theory is the following one.

\begin{thm}
Let $\mathcal{C}$ be an abelian category endowed with
a functor
$F\colon\mathcal{C}\to\Vect$ whose values are finite
dimensional vector spaces and which is both exact and
faithful. The category $\mathcal{C}$ is equivalent
to the category $\Com\EndCheck(F)$ of finite dimensional
$\EndCheck(F)$\nobreakdash-\hspace{0pt}comodules. Moreover, if
$\mathcal{C}$ is equivalent as a Tannakian category to $\Com C$
for another coalgebra $C$, then $C$ is canonically isomorphic 
to $\EndCheck(F)$.
\end{thm}

\section{Tannaka Reconstruction for Crossed Structures}
We introduce now the analog of Tannaka Theory in the context
of crossed structures.

Let $\mathcal{C}$ be an autonomous abelian 
crossed $G$\nobreakdash-\hspace{0pt}category (again we suppose $\mathcal{C}$
strict, but the results can obviously be generalized to the case of a category
equivalent to a strict one).
A \textit{fiber crossed  $G$\nobreakdash-\hspace{0pt}functor}
is a couple $\bigl(F,\{\dot{\varphi}_{\beta}\}_{\beta\in G}\bigr)$
such that $F\colon\mathcal{C}\to\Vect$ is an autonomous tensor functor
(i.e., a tensor functor with preserves dualities) whose values 
are finite dimensional vector spaces and which is 
faithful and exact, while 
$\dot{\varphi}_{\beta}\colon F\to F\circ\Phi_{\beta}$ 
is a natural isomorphism satisfying the following conditions.
\begin{itemize}
\item $\dot{\varphi}_{\beta_{1}}\circ\dot{\varphi}_{\beta_{2}}=
      \dot{\varphi}_{\beta_{1}\beta_{2}}$ for all
      $\beta_{1},\beta_{2}\in G$ and $\dot{\varphi}_{1_{G}}=\Id$.
\item If we set $F_{\alpha}=F\vert_{\mathcal{C}_{\alpha}}$, then
      $\dot{\varphi}_{\beta}$ induces a natural isomorphism
      $F_{\alpha}\cong F_{\beta\alpha\beta^{-1}}\circ\Phi_{\beta}$.
\end{itemize}
We say that $\mathcal{C}$ is a 
\textit{Tannakian crossed $G$\nobreakdash-\hspace{0pt}category} 
if it is endowed with a fiber functor. We say that two 
Tannakian crossed $G$\nobreakdash-\hspace{0pt}categories $\mathcal{C}$ 
with fiber functor $X$ and $\mathcal{D}$ with fibre functor $Y$
are \textit{equivalent} [resp.\@ \textit{isomorphic\/}] if there is an
equivalence [resp. an isomorphism] of crossed 
$G$\nobreakdash-\hspace{0pt}categories $F\colon\mathcal{C}\to\mathcal{D}$
such that $X=Y\circ F$.

\begin{thm}
Let $\mathcal{C}=
\bigl(\mathcal{C},F,\{\dot{\varphi}_{\beta}\}_{\beta\in G}\bigr)$ 
be a Tannakian crossed $G$\nobreakdash-\hspace{0pt}category.
There exists a crossed Hopf
$G$\nobreakdash-\hspace{0pt}algebra $H=H(\mathcal{C})$,
unique up to isomorphism, such that
$\mathcal{C}\cong\Com H$ as Tannakian crossed
$G$\nobreakdash-\hspace{0pt}categories.
\end{thm}

\begin{sketch}
For all $\alpha\in G$, let us set $H_{\alpha}=\EndCheck(F_{\alpha})$.
The functor $F_{\alpha}$ is obviously exact and faithful, so we have
$\mathcal{C}_{\alpha}\cong\Com H_{\alpha}$. By Lemma~\ref{l:trivial},
the conjugation $\Phi$ of $\mathcal{C}$ and the $\dot{\varphi}$
give rise to a family of isomorphism 
$\varphi_{\beta}\colon H_{\alpha}\to H_{\beta\alpha\beta^{-1}}$
with the obvious property that $\varphi_{1_{G}}=\Id$ and
$\varphi_{\beta_{1}\beta_{2}}=\varphi_{\beta_{1}}\circ\varphi_{\beta_{2}}$.

Suppose, for simplicity, that $\mathcal{C}$ is strict. Then, by 
Lemma~\ref{l:noja}, with $X=F_{\alpha}$ and $Y=F_{\beta}$, and observing
that $F_{\alpha}(U)\otimes F_{\beta}(V)= F_{\alpha\beta}(U\otimes V)$
for all objects $U$ and $V$ in $\mathcal{C}$, we get a morphism
(which is actually also an isomorphism)
\begin{equation*}
m_{\alpha,\beta}\colon H_{\alpha}\otimes 
H_{\beta}\to H_{\alpha\beta}\text{.}
\end{equation*}
It is easy to prove that in that way we obtained a multiplication
for $H$ with unit $[1]$, being $1$ the element of $F_{1_{G}}(\tunit)$
corresponding to $1\in\C$ under the isomorphism 
$F_{1_{G}}(\tunit)\cong\C$. If $\mathcal{C}$ is not strict,
then one can define the product of $[h]$ and $[k]$ (where
$h\in\End\bigl(F(U)\bigr)$ and  $k\in\End\bigl(F(V)\bigr)\/$)
as the class of the map
\begin{equation*}
 F(U\otimes V)\cong F(U)\otimes F(V)\xrightarrow{h\otimes k}
 F(U)\otimes F(V)\cong F(U\otimes V)
\end{equation*}
(notice that exactly the way used by Joyal and Street~\cite{{JS-Tannaka}}
to endow the coalgebra $H=\bigoplus_{\alpha\in G}H_{\alpha}$ of a structure
of bialgebra).

Now, let $U$ be an object in $\mathcal{C}_{\alpha}$ and let $U^{\ast}$ be
its left dual. Since tensor functors preserve the dual pairings, for all
$h\in\End\bigl(F_{\alpha}(U)\bigr)$ there is a transposed endomorphism
${}^{t}h\in\End\bigl(F_{\alpha^{-1}}(U^{\ast})\bigr)$. We set
\begin{equation*}
 s_{\alpha}\bigl([h]\bigr)=\bigl[{}^{t}h\bigr]\text{.}
\end{equation*}
The proof that $s$ is an antipode is then an adaptation of that 
in~\cite{JS-Tannaka} (Proposition~5, page~468) and is omitted here.

Finally, by using Lemma~\ref{l:trivial}, we get a family of isomorphisms
$\varphi_{\beta}\colon H_{\alpha}\to H_{\beta\alpha\beta^{1}}$. To check
that $\varphi$ satisfies our axioms is routine.
\end{sketch}

\begin{cor} If $\mathcal{C}$ is braided, then $H$ is coquasitriangular
and $\mathcal{C}\cong\Com H$ as braided crossed
$G$\nobreakdash-\hspace{0pt}categories. If $\mathcal{C}$ is ribbon, 
then $H$ is coribbon
and $\mathcal{C}\cong\Com H$ as ribbon crossed
$G$\nobreakdash-\hspace{0pt}categories.
\end{cor}

\begin{proof}
Suppose that $\mathcal{C}$ is braided. The equivalence
$E\colon\Com H\to\mathcal{C}$ induces a braiding in $\Com H$
by setting, for all objects $M$ and $N$ in $\Com H$,
$c_{M,N}=E^{-1}(c_{E(M),E(N)})$. Thus, by Theorem~\ref{thm:Beethoven},
$H$ is coquasitriangular. The proof for the ribbon is similar.
\end{proof}

\section{Quantum Co-double Construction}
In~\cite{Zunino-1} the author generalized the center construction
for tensor categories~\cite{JS} to the case of crossed group categories.
Stating from any crossed group category $\mathcal{C}$, the center 
construction provides a braided tensor category $\mathcal{Z}(\mathcal{C})$.
In similar way, by means of the generalization in~\cite{Zunino-1} of the
construction in~\cite{KasTur} and~\cite{Street-double}, starting from
a braided crossed $G$\nobreakdash-\hspace{0pt}category $\mathcal{D}$,
one can construct
a ribbon crossed $G$\nobreakdash-\hspace{0pt}category
$\theta(\mathcal{D})$. By means
of the Tannaka theory for crossed $G$\nobreakdash-\hspace{0pt}categories,
we can use
these constructions to recover in the crossed case the results about the
quantum co-double of a Hopf algebra discussed in~\cite{Street-double}.
Let us recall the definition of both the center $\mathcal{Z}(\cdot)$
and the ribbon extension $\theta(\cdot)$ as in~\cite{Zunino-1}. 

\subsection{The Center} Let $\mathcal{C}$ be a 
crossed $G$\nobreakdash-\hspace{0pt}category. The component 
$\mathcal{Z}_{\alpha}(\mathcal{C})=\bigl(\mathcal{Z}(\mathcal{C})\bigr)_{\alpha}$ 
of $\mathcal{Z}(\mathcal{C})$ is the category whose
objects are couples $(U,\mathfrak{c}_{\_})$, where $U$ is an object in 
$\mathcal{C}_{\alpha}$ and $\mathfrak{c}_{\_}$ is a natural isomorphism
of endomorphic functors in $\mathcal{C}$
\begin{equation*}
\mathfrak{c}_{\_}\colon \Bigl(U\otimes\_\Bigr)\to\lidx{U}{\_}\otimes U
\end{equation*}
such that, for all $V,W\in\mathcal{C}$,
\begin{equation*}
\mathfrak{c}_{V\otimes W}=\biggl(\Bigl(\lidx{U}{V}\Bigr)
\otimes\mathfrak{c}_{W}\biggr)
\otimes(\mathfrak{c}_{V}\otimes Y)\text{.}
\end{equation*}
If both $(U,\mathfrak{c}_{\_})$ and 
$(V,\mathfrak{d}_{\_})$ are objects in $\mathcal{Z}_{\alpha}(\mathcal{C})$,
then an arrow $f\colon (U,\mathfrak{c}_{\_})\to (V,\mathfrak{d}_{\_})$ 
is a map
$f\colon U\to V$ such that 
\begin{equation*}
 \biggl(\Bigl(\lidx{W}{V}\Bigr)\otimes 
 f\biggr)\circ\mathfrak{c}_{W}=\mathfrak{d}_{W}\circ
(f\otimes W)\text{,}
\end{equation*}
for all $W\in\mathcal{C}$. 
The tensor product of $(U,\mathfrak{c})$ in $\mathcal{Z}_{\alpha}(\mathcal{C})$
and $(U',\mathfrak{c}')$ in $\mathcal{Z}_{\alpha'}(\mathcal{C})$
is obtained by setting
\begin{equation*}
(U,\mathfrak{c})\otimes (U',\mathfrak{c'})=
\bigl(U\otimes U', \mathfrak{c}\overline{\otimes}\mathfrak{c'}\bigr)\text{,}
\end{equation*}
where, for all $W\in\mathcal{C}$,
\begin{equation*}
 (\mathfrak{c}\overline{\otimes}\mathfrak{c}')_{W}=
 (\mathfrak{c}_{\Bigl(\lidx{V}{W}\Bigr)}\otimes V)
 \circ(U\otimes\mathfrak{c}'_{W})\text{.}
\end{equation*}
Finally $\Phi_{\beta}(U,\mathfrak{c}_{\_})=
\Bigl(\Phi_{\beta}(U),\mathfrak{c}^{\beta}_{\_}\Big)$,
where, for all $V\in\mathcal{C}$, 
\begin{equation*}
\mathfrak{c}^{\beta}_{V}=
\Phi_{\beta}(\mathfrak{c}_{\varphi^{-1}_{\beta}(X)})\text{.}
\end{equation*}
The crossed $G$\nobreakdash-\hspace{0pt}category $\mathcal{Z}(\mathcal{C})$
is braided by setting
\begin{equation*}
 c_{(U,\mathfrak{c}_{\_}),(U',\mathfrak{c}'_{\_})}=
 \mathfrak{c}_{U'}\text{.}
\end{equation*}

\begin{lemma}\label{l:final}
If $\mathcal{C}$ is abelian, then $\mathcal{Z}(\mathcal{C})$
is abelian and the forgetful functor
$\map{Z}{\mathcal{Z}(\mathcal{C})}{\mathcal{C}}%
 {(U,\mathfrak{c}_{\_})}{U}$ is exact.
\end{lemma}

The following proof is adapted from the one for the center
of a tensor category given by Street~\cite{Street-double}.

\begin{proof}
 For all $\alpha\in G$, let $0_{\alpha}$ be the object $0$ of
 $\mathcal{C}_{\alpha}$. For all objects $X$ in $\mathcal{C}_{\beta}$,
 we have
 \begin{equation*}
  0_{\alpha}\otimes X = 0_{\alpha\beta}
  =\Bigl(\lidx{\alpha}{X}\Bigr)\otimes 0_{\beta}\text{.}
 \end{equation*}
So, if $0_{\alpha}[X]$ is the identity of $0_{\alpha\beta}$, then we
get an object $\bigl(0_{\alpha},0_{\alpha}[\cdot]\bigr)$ in
$\mathcal{Z}(\mathcal{C})$.

Let both $(U,\mathfrak{c}_{\_})$ and $(U',\mathfrak{c}'_{\_})$ be
objects in $\mathcal{Z}(\mathcal{C})$ and let $X$ be an 
object in $\mathcal{C}$. By setting
\begin{equation*}
(\mathfrak{c}\oplus\mathfrak{c}')_{X}\colon
(U\oplus U')\otimes X\cong (U\otimes X)\oplus (U'\otimes X)
\xrightarrow{\mathfrak{c}_{X}\oplus\mathfrak{c}'_{X}}
\biggl(\Bigl(\lidx{\alpha}{X}\Bigr)\otimes U\biggr)\oplus
\biggl(\Bigl(\lidx{\alpha}{X}\Bigr)\otimes U'\biggr)\cong
\Bigl(\lidx{\alpha}{X}\Bigr)\otimes (U\oplus U')\text{,}
\end{equation*}
we get an object $(U,\mathfrak{c}_{\_})\oplus (U',\mathfrak{c}'_{\_})=
\bigl(U\oplus U',(\mathfrak{c}\oplus\mathfrak{c}')_{\_}\bigr)$ in
$\mathcal{C}_{\alpha}$.

Let $f\colon (U,\mathfrak{c}_{\_})\to (U',\mathfrak{c}'_{\_})$
be an arrow in $\mathcal{Z}_{\alpha}(\mathcal{C})$ and let
$\krn(f)\colon\Krn(f)\to U$ and 
$\cok(f)\colon U'\to\Cok(f)$ be the kernel and the cokernel of $f$
in $\mathcal{C}_{\alpha}$. Since the following diagram in exact
\begin{equation*}
\vcenter{\xymatrix@C=8ex{
 \ \ \ \protect\Krn(f)\otimes X\ \ \   
 \ar@{>->}[r]^{\protect\krn(f)\otimes X}
 \ar@{.>}[d]_(.4){\protect\mathfrak{K}[f]_{X}\ }
 &
 U\otimes X 
 \ar[r]^{f\otimes X}
 \ar[d]|(.4){\mathfrak{c}_{X}}
 &
 U'\otimes X
 \ar[d]|(.4){\mathfrak{c}'_{X}}
 \ar@{->>}[r]^{\protect\cok(f)\otimes X}
 &
 \protect\Cok(f)\otimes X 
 \ar@{.>}[d]^(.4){\ \protect\mathfrak{C}[f]_{X}} \\
 \ \ \ \protect\Bigl(\protect\lidx{\alpha}{X}\protect\Bigr)\otimes
 \protect\Krn(f)\ \ \ 
 \ar@{>->}[r]_{\ \ \protect\bigl(\protect\lidx{\alpha}{X}\protect\bigr)
         \otimes\protect\krn(f)}
 &
 \protect\Bigl(\protect\lidx{\alpha}{X}\protect\Bigr)\otimes U
 \ar[r]_{\protect\bigl(\protect\lidx{\alpha}{X}\bigr)\otimes f} &
 \protect\Bigl(\protect\lidx{\alpha}{X}\protect\Bigr)\otimes U'
 \ar@{->>}[r]_{\protect\bigl(\protect\lidx{\alpha}{X}\protect\bigr)
  \otimes\protect\cok(f)} &
 \protect\Bigl(\protect\lidx{\alpha}{X}\protect\Bigr)\otimes
 \protect\Cok(f)
 }}\text{,}
\end{equation*}
there are unique $\mathfrak{K}^{f}_{X}$ and $\mathfrak{C}^{f}_{X}$
making it commutative. One can prove that 
$\krn(f)\colon\bigl(\Krn(f)\otimes U,\mathfrak{K}[f]_{\_}\bigr)\to
(U,\mathfrak{c}_{\_})$ and $\cok(f)\colon (U',\mathfrak{c}'_{\_})\to
\bigl(\Cok(f),\mathfrak{C}[f]_{\_}\bigr)$ are the kernel and the 
cokernel of $f$ in $\mathcal{Z}_{\alpha}(\mathcal{C})$.
The rest follows easily.
\end{proof}

\subsection{The Ribbon Extension} Let $\mathcal{D}$ be a
braided crossed $G$\nobreakdash-\hspace{0pt}category. To define
$\theta(\mathcal{D})$ let us firstly define the balanced 
crossed $G$\nobreakdash-\hspace{0pt}category $\mathcal{D}^{Z}$.
The component 
$\mathcal{D}^{Z}_{\alpha}=\Bigl(\mathcal{D}^{Z}\Bigr)_{\alpha}$ 
of $\mathcal{D}^{Z}$ is the category whose
objects are couples $(U,t)$, where $U$ is an object in 
$\mathcal{D}_{\alpha}$ and $t\colon U\to\lidx{U}{U}$
is an isomorphism in $\mathcal{D}$.
An arrow $f\colon (U,t)\to(U',t')$ in $\mathcal{D}^{Z}$ is a map
$f\colon U\to U'$ is $\mathcal{D}_{\alpha}$ such that
\begin{equation*}
\Bigl(\lidx{U}{f}\Bigr)\circ t= t'\circ f\text{.}
\end{equation*}
We set
\begin{equation*}
(U,t)\otimes (U',t')
=(U\otimes U',t\hat{\otimes} t')\text{,}
\end{equation*}
where
\begin{equation*}
t\hat{\otimes}t'=
c_{\Bigl(\leftidx{^{U\otimes U'}}{U}{'}\Bigr),\Bigl(\lidx{U}{U}\Bigr)}\circ
c_{\Bigl(\lidx{U}{U}\Bigr),\Bigl(\leftidx{^{U'}}{U}{'}\Bigr)}\circ
(t\otimes t')\text{.}
\end{equation*}
Finally, we set
\begin{equation*}
 \Phi_{\beta}(U, t)=
 \bigl(\Phi_{\beta}(U),\Phi_{\beta}(t)\bigr)\text{.}
\end{equation*}
The braiding of $\mathcal{D}$ induces a braiding in $\mathcal{D}^{Z}$
and $\mathcal{D}^{Z}$ is balanced with the twist
\begin{equation*}
 \theta_{(U, t)}=t\text{.}
\end{equation*}

The following lemma follows by the proof of
Proposition~4 in~\cite{Street-double}.

\begin{lemma}
If $\mathcal{D}$ is abelian, then $\mathcal{D}^{Z}$
is abelian and the forgetful functor
$\map{Z_{1}}{\mathcal{D}^{Z}}{\mathcal{D}}{(U,t)}{U}$
is exact.
\end{lemma}

The crossed
$G$\nobreakdash-\hspace{0pt}category 
$\mathcal{D}^{Z}$ is not necessarily ribbon. However,
let $\mathcal{E}$ balanced crossed
$G$\nobreakdash-\hspace{0pt}category and let $\mathcal{N}(\mathcal{E})$
the full subcategory of $\mathcal{E}$
of objects $U$ such that there exists a left dual $U^{\ast}$ 
\begin{itemize}
\item compatible with $\Phi$, i.e. such
that eq.\@~\eqref{e:good} is satisfied,
\item compatible with $\theta$, i.e. such that eq.\@~\eqref{e:tortility}
is satisfied, and
\item we have
\begin{equation*}
\theta^{-2}_{U}=\omega_{(b_{U},d_{U})}\text{,}
\end{equation*}
where, by definition,
$\theta^{2}_{U}=\bigl(\leftidx{^{U}}{\theta}{_{U}}\bigr)\circ\theta_{U}$,
$\theta^{-2}_{U}=(\theta^{2}_{U})^{-1}$, and
\begin{equation*}
 \omega_{(b_{U},d_{U})}=(d_{\leftidx{^{U\otimes U}}{U}{^{\ast}}}\otimes U)\circ
 \biggl(\Bigl(\leftidx{^{U\otimes U}}{U}{^{\ast}}\Bigr)\otimes
 \tilde{c}_{\lidx{U}{U}, \lidx{U\otimes U}{U}}\biggr)\circ
 \Bigl((c_{\,\lidx{U}{U},\leftidx{^{U}}{\! U}{^{\ast}}}\circ
  b_{\lidx{U}{U}})\otimes\lidx{U\otimes U}{U}\Bigr)\text{,}
\end{equation*}
where, for all objects $X$ in $\mathcal{C}_{\alpha}$ and
$Y$ in $\mathcal{C}_{\beta}$, we set
 $\tilde{c}_{X,Y}=(c_{Y,\lidx{\beta^{-1}}{X}})^{-1}$.
\end{itemize}
The category 
$\mathcal{N}(\mathcal{E})$ is a ribbon crossed
$G$\nobreakdash-\hspace{0pt}category. Thus we set
\begin{equation*}
 \theta(\mathcal{D})=\mathcal{N}\Bigl(\mathcal{D}^{Z}\Bigr)\text{.}
\end{equation*}

\begin{lemma}
If $\mathcal{E}$ is abelian, then $\mathcal{N}(\mathcal{E})$
is abelian and the inclusion 
$N\colon\mathcal{N}(\mathcal{E})\to\mathcal{E}$ is exact. 
So, when $\mathcal{E}=\mathcal{D}^{Z}$ we get an exact functor
$Z_{2}=Z_{1}\circ N\colon\theta(\mathcal{D})\to\mathcal{D}$.
\end{lemma}

The following proof is again an adaptation of that given in 
Street~\cite{Street-double} in the case of a tensor category.

\begin{proof}
We only need to show that $\mathcal{N}(\mathcal{E})$
is closed under finite limits and finite colimits.
Since $(U\oplus V)^{\ast}= U^{\ast}\oplus V^{\ast}$
and $\theta_{(U\oplus V)^{\ast}}=\theta^{\ast}_{U\otimes V}$, 
the category $\mathcal{N}(\mathcal{E})$ is closed under
direct sums. Now, if $f\colon U\to V$ is a map in an
autonomous category and both $U$ and $V$ has left
duals, then the cokernel object 
$\Cok(f^{\ast})$ of $f^{\ast}$ is a left dual of the 
kernel object $\Krn(f)$ of $f$ while the kernel object
$\Krn(f^{\ast})$ of $f^{\ast}$ is a right dual of the
cokernel object $\Cok(f)$ of $f$. Further, if $f$ is a
map in an abelian crossed $G$\nobreakdash-\hspace{0pt}category
and eq.\@~\eqref{e:good} holds for both $U^{\ast}$ and $V^{\ast}$
then, it holds also for both $\Cok(f^{\ast})$ and $\Krn(f^{\ast})$,
since the $\Phi_{\alpha}$ are exact. If we are in a braided 
crossed $G$\nobreakdash-\hspace{0pt}category, then a right dual
is also a left dual, and so $\Krn(f^{\ast})$ is a left dual
of $\Cok(f)$. So, if $f$ is a map in $\mathcal{N}(\mathcal{E})$,
both $\Krn(f)$ and $\Cok(f)$ have left dual in $\mathcal{E}$.
By using the naturality of the twist, one can prove that both
$\Krn(f)$ and $\Cok(f)$ also satisfy eq.\@~\eqref{e:tortility}.
\end{proof}

\subsection{The Quantum Co-double}

\begin{thm}
 Let $H$ be a crossed Hopf
 $G$\nobreakdash-\hspace{0pt}algebra. There exists a 
 coquasitriangular crossed Hopf 
 $G$\nobreakdash-\hspace{0pt}algebra $D^{\ast}(H)$
 \textup{(\/}unique up to isomorphism\textup{)}
 such that 
 \begin{equation*}
  \Com\bigl(D^{\ast}(H)\bigr)=\mathcal{Z}(\Com H)
 \end{equation*} 
 as braided Tannakian 
 crossed $G$\nobreakdash-\hspace{0pt}categories.
 Let $H'$ be a coquasitriangular crossed Hopf
 $G$\nobreakdash-\hspace{0pt}algebra. There exists a 
 coribbon crossed Hopf 
 $G$\nobreakdash-\hspace{0pt}algebra $R(H')$ 
 \textup{(\/}again unique up to isomorphism\textup{)}
 such that 
 \begin{equation*}
  \Com\bigl(R(H')\bigr)=\theta(\Com H')
 \end{equation*}
 as ribbon Tannakian 
 crossed $G$\nobreakdash-\hspace{0pt}categories. 
 In particular, for $H'=D^{\ast}(H)$ we get
 \begin{equation*}
  \Com\Bigl(R\bigl(D^{\ast}(H)\bigr)\Bigr)=
  \theta\bigl(\mathcal{Z}(\Com H)\bigr)\text{.}
 \end{equation*}
\end{thm}

\begin{proof}
Let $\mathcal{C}=\Com H$ and let $\vert\cdot\vert\colon \Com H\to\Vect$
and $Z\colon\mathcal{Z}(\mathcal{C})\to\mathcal{C}$ be
the trivial forgetful functors. By Lemma~\ref{l:final}, the functor 
$Z\circ\vert\cdot\vert\colon\mathcal{Z}(\Com H)\to\Vect$ 
is a fiber crossed  $G$\nobreakdash-\hspace{0pt}functor.
Thus, we are allowed to apply Tannaka theory
getting the Hopf $G$\nobreakdash-\hspace{0pt}algebra $D(H)$ such that 
$\mathcal{Z}(\Com H)\cong \Com\bigl(D(H)\bigr)$

The proof for the ribbon structures is similar and left to the reader.
\end{proof}

We recall that, given a crossed Hopf $G$\nobreakdash-\hspace{0pt}coalgebra $K$ of
finite type, one can construct a quasitriangular crossed Hopf 
$G$\nobreakdash-\hspace{0pt}coalgebra $\overline{D}(K)$ such that 
$\mathcal{Z}\bigl(\Mod K\bigr)=\Mod\bigl(\overline{D}(K)\bigr)$, being $\Mod K$
the category of finite dimensional $K$\nobreakdash-\hspace{0pt}modules, 
see~\cite{Zunino-1,Zunino-2}. 
In particular,
$\overline{D}(K)$ is of finite type if and only if $H$ is totally finite
(i.e. $\sum_{\alpha\in G}\dim K_{\alpha}<+\infty\/$). In that case, both the dual of
$K$ and the dual of $\overline{D}(K)$ are crossed Hopf $G$\nobreakdash-\hspace{0pt}algebras.
It is easy to prove that we have an equivalence of crossed
$G$\nobreakdash-\hspace{0pt}categories $\Mod K\cong \Com K^{\ast}$. We deduce the
following chain of equivalence of braided crossed group categories
\begin{equation}\label{e:congu}
\Com\Bigl(\bigl(\overline{D}(K)\bigr)^{\ast}\Bigr)=
\Mod\bigl(\overline{D}(K)\bigr)\cong
\mathcal{Z}(\Mod K)\cong\mathcal{Z}\bigl(\Com(K^{\ast})\bigr)\cong
\Com\bigl(D^{\ast}(K^{\ast})\bigr)\text{.}
\end{equation}

\begin{cor}\label{cor:1}
Let $H$ be a totally finite crossed Hopf $G$\nobreakdash-\hspace{0pt}algebra. The
coquasitriangular crossed Hopf $G$\nobreakdash-\hspace{0pt}algebras $D^{\ast}(H)$
and $\bigl(\overline{D}(H^{\ast})\bigr)^{\ast}$ are isomorphic.
\end{cor}

\begin{proof}
By setting $K=H^{\ast}$ in~\eqref{e:congu} we get an equivalence of
braided $G$\nobreakdash-\hspace{0pt}categories 
\begin{equation*}
\Com\Bigl(\bigl(\overline{D}(H^{\ast})\bigr)^{\ast}\Bigr)=
\Com\bigl(D^{\ast}(H)\bigr)\text{.}
\end{equation*} By Tannaka theory, the two
crossed Hopf $G$\nobreakdash-\hspace{0pt}algebras must be isomorphic.
\end{proof}

Starting from a braided crossed Hopf $G$\nobreakdash-\hspace{0pt}coalgebra $L$ 
of finite type, there exists a ribbon crossed Hopf $G$\nobreakdash-\hspace{0pt}coalgebra
$RT(L)$ such that $\Mod\bigl(RT(L)\bigr)\cong\theta(\Mod L)$, 
see~\cite{Zunino-1,Zunino-2}. By means of a chain of equivalence of 
ribbon $G$\nobreakdash-\hspace{0pt}categories similar to the previous one we get
the following result.

\begin{cor}\label{cor:2}
Let $H'$ be a Hopf $G$\nobreakdash-\hspace{0pt}algebra of finite type. The
coribbon crossed Hopf $G$\nobreakdash-\hspace{0pt}algebras $R(H')$ and 
$\bigl(RT(H^{\ast})\bigr)^{\ast}$ are isomorphic.
\end{cor}

By combining the two previous corollaries, we finally get the
following result.

\begin{cor}
Let $H$ be a totally finite crossed Hopf $G$\nobreakdash-\hspace{0pt}algebra. The
coribbon crossed Hopf $G$\nobreakdash-\hspace{0pt}algebras 
$R\bigl(D^{\ast}(H)\bigr)$
and $\Bigl(RT\bigl(\overline{D}(H^{\ast})\bigr)\Bigr)^{\ast}$ are isomorphic.
\end{cor}

The quantum double of a totally finite semisimple crossed Hopf 
$G$\nobreakdash-\hspace{0pt}coalgebra is modular with 
$\theta_{\alpha}=1_{\alpha}$~\cite{Zunino-1}. One can easily prove 
the following result.

\begin{cor}
The quantum co-double $D^{\ast}(H)$ of a totally finite crossed Hopf 
$G$\nobreakdash-\hspace{0pt}algebra is modular with
$\tau_{\alpha}=\varepsilon_{\alpha}$.
\end{cor}

\end{document}